\documentclass[12pt]{amsart}
\input{amssymb.sty}
\title[Some prime factorization results 
for type $\mathrm{II}_1$ factors]
{Some prime factorization results\\ 
for type $\mathrm{II}_1$ factors}
\author[N. Ozawa]{Narutaka OZAWA$^*$}
\address{Department of Mathematical Sciences,
University of Tokyo, Komaba, 153-8914}
\email{narutaka@ms.u-tokyo.ac.jp}
\thanks{${}^*$ Supported by the JSPS 
Postdoctoral Fellowships for Research Abroad}
\author[S. Popa]{Sorin Popa$^{**}$}
\address{Department of Mathematics, 
UCLA, Los Angeles, CA90095}
\email{popa@math.ucla.edu}
\thanks{${}^{**}$ Supported by NSF-Grant 0100883}
\date{Feburuary 27, 2003}

\newtheorem{thm}{Theorem}
\newtheorem{prop}[thm]{Proposition}
\newtheorem{cor}[thm]{Corollary}
\newtheorem{lem}[thm]{Lemma}
\newcommand{\G}{\Gamma} 
\newcommand{\e}{\varepsilon} 
\newcommand{\N}{{\mathbb N}} 
\newcommand{\C}{{\mathbb C}} 
\newcommand{\B}{{\mathbb B}} 
\newcommand{\K}{{\mathbb K}} 
\newcommand{\M}{{\mathbb M}} 
\newcommand{\F}{{\mathbb F}} 
\newcommand{\R}{{\mathbb R}}
\newcommand{\CC}{{\mathcal C}}
 
\newcommand{\FF}{{\mathcal F}} 
\newcommand{\LG}{{\mathcal L}\G} 
\newcommand{\LL}{{\mathcal L}} 
\newcommand{\MM}{{\mathcal M}} 
\newcommand{\NN}{{\mathcal N}} 
\newcommand{\PP}{{\mathcal P}} 
\newcommand{\QQ}{{\mathcal Q}} 
\newcommand{\RR}{{\mathcal R}} 
\newcommand{\ZZ}{{\mathcal Z}} 
\newcommand{\A}{{\mathcal A}} 
\newcommand{\BB}{{\mathcal B}} 
\newcommand{\U}{{\mathcal U}} 
\newcommand{\Proj}{\mathrm{Proj}}
 
\newcommand{\hh}{{\mathcal H}} 
 
\newcommand{\eg}{\ell_2\G} 
\newcommand{\id}{\mathrm{id}} 
\newcommand{\Out}{\mathop{\mathrm{Out}}} 
\newcommand{\Mod}{\mathop{\mathrm{Mod}}} 
\newcommand{\MOD}{\mathop{\mathrm{\bf Mod}}} 
\newcommand{\Ad}{\mathop{\mathrm{Ad}}} 
\newcommand{\om}{\otimes_{\min}} 
\newcommand{\conv}{\overline{\mathrm{co}}^w} 
 
\newcommand{\clg}{C^*_\lambda\G}
\newcommand{\crg}{C^*_\rho\G}
\newcommand{\clrg}{C^*_{\lambda,\rho}\G}
\newcommand{\ov}{\overline{\otimes}}

\begin{document}
\begin{abstract} We prove several unique prime 
factorization results for tensor products of 
type II$_1$ factors coming from groups that can be 
realized either as subgroups of 
hyperbolic groups or as discrete subgroups of connected 
Lie groups of real rank 1. In particular, we show 
that if $\RR\ov\LL\F_{r_1}\ov\cdots\ov\LL\F_{r_m}$ 
is isomorphic to a subfactor in 
$\RR\ov\LL\F_{s_1}\ov\cdots\ov\LL\F_{s_n}$, 
for some $2\leq r_i, s_j \leq \infty$,  
then $m\le n$. 
\end{abstract}
\maketitle
\section{Introduction}
Type $\mathrm{II}_1$ factors coming from non-amenable 
hyperbolic groups were shown 
to be prime in \cite{ozawa}, while a uniqueness result,  
modulo unitary conjugacy, for certain Cartan subalgebras 
in type $\mathrm{II}_1$ factors was established in \cite{popab}.  
In this paper we will combine ideas and techniques from \cite{ozawa} 
and \cite{popab} to prove a variety of 
unique prime factorization results, modulo unitary conjugacy, 
for tensor products of 
type $\mathrm{II}_1$ factors coming from hyperbolic groups. 
To state the results we need some notations.   

Thus, we denote by $\CC$ 
the class of all countable non-amenable ICC 
groups which are either (subgroups of) hyperbolic groups 
or discrete subgroups of connected simple Lie groups 
of rank one. We denote by 
$\RR$ the hyperfinite type $\mathrm{II}_1$ factor. 

If $\MM$ is a type $\mathrm{II}_1$factor 
and $t > 0$ 
then we denote by $\MM^t$ the amplification 
of $\MM$ by $t$, i.e., $\MM^t = p\M_{n\times n}(\MM)p$ 
for some $n \geq t$ and 
$p \in\M_{n\times n}(\MM)$ a projection 
of trace $\tau(p)=t/n$. Thus, $\MM^t$ is defined up to 
isomorphism, in fact up to unitary conjugacy in some 
``large matrix'' algebra over $\MM$. 

Similarly, if a type $\mathrm{II}_1$ factor $\MM$ is decomposed as 
$\MM=\MM_1 \ov \MM_2$ for some type $\mathrm{II}_1$ factors 
$\MM_1, \MM_2$ 
and $t > 0$ then, modulo unitary conjugacy, there exists a 
unique decomposition 
$\MM=\MM_1^t \ov \MM_2^{1/t}$, 
such that $p_1\MM_1p_1 \vee p_2\MM_2p_2$ and 
$q_1\MM_1^tq_1 \vee q_2\MM_2^{1/t}q_2$ are unitary conjugate in $\MM$ 
for any projections 
$p_i \in \Proj(\MM_i)$, $i=1,2$ and $q_1 \in \Proj(\MM_1^t)$, 
$q_2 \in \Proj(\MM_2^{1/t})$, 
with $\tau(p_1)/\tau(q_1)=\tau(q_2)/\tau(p_2)=t.$

With these notations in mind we have: 

\begin{thm}\label{split1}
Let $\G_1,\cdots,\G_n\in\CC$ and 
$\G=\G_1\times\cdots\times\G_n$. 
Assume that $\LG = \NN_1\ov\NN_2$ 
for some type $\mathrm{II}_1$ factors $\NN_1$ and $\NN_2$. 
Then, there exist $t>0$ and 
a partition $I_1\cup I_2=\{1,\ldots,n\}$ such that 
$\NN_1^t = \overline{\bigotimes}_{k\in I_1}\LG_k$ 
and 
$\NN_2^{1/t}=\overline{\bigotimes}_{k\in I_2}\LG_k$ 
modulo conjugacy by a unitary element in $\LG$. 
\end{thm}

An analogous result holds true for McDuff factors: 

\begin{thm}\label{split2}
Let $\G_1,\cdots,\G_n\in\CC$ and 
$\G=\G_1\times\cdots\times\G_n$. 
Assume that $\RR\ov\LG = \NN_1\ov\NN_2$ 
for some McDuff type $\mathrm{II}_1$ factors $\NN_1$ and $\NN_2$.  
Then, there exist a partition $I_1\cup I_2=\{1,\ldots,n\}$ 
and a decomposition $\RR = \RR_1 \ov \RR_2$ such that 
$\NN_1 = \RR_1\ov\overline{\bigotimes}_{k\in I_1}\LG_k$ 
and 
$\NN_2 = \RR_2\ov\overline{\bigotimes}_{k\in I_1}\LG_k$ 
modulo unitary conjugacy. 
Moreover, $I_1$ (resp.\ $I_2$) is empty if and only if 
$\NN_1$ (resp.\ $\NN_2$) is injective. 
\end{thm}

The proofs of both theorems follow the same strategy: 
A careful use of the ideas 
and techniques in \cite{ozawa} shows 
that if $\NN=\NN_1\ov \NN_2$ is a subfactor 
of the group von Neumann algebra $\MM=\MM_1 \ov \MM_2$, with 
$\MM_i=\LG_i$ and $\G_i$ 
hyperbolic groups, then either $\NN_1$ or $\NN_2$ 
must ``avoid killing'' the compacts $K (L^2(\MM_i)) \otimes \C 1$. 
This provides a non-zero compact (relative to $\MM_i$) element in the  
commutant of say $N_1$, thus a Hilbert 
$\NN_1-\MM_i$ bimodule 
$\hh \subset L^2(\MM)$ with dim$\hh_{\MM_i}< \infty$ (cf. \cite{popab}). 
Using techniques from \cite{popab} this is shown to imply 
that $\NN_1^t$ is unitary conjugate to $\MM_i$ for some $t > 0$. 

The following four corollaries are immediate consequences 
of the above theorems. 

\begin{cor}
Let $\G_1,\cdots,\G_n\in\CC$. 
Assume 
\[
\NN_1\ov\cdots\ov\NN_m=\LG_1\ov\cdots\ov\LG_n, 
\]
for some $m\geq n$ and type $\mathrm{II}_1$ factors 
$\NN_1,\cdots,\NN_m$. 
Then $m=n$ and there exist 
$t_1, t_2, ..., t_n > 0$ 
with $t_1t_2\cdots t_n = 1$ such that after permutation of indices 
and unitary conjugacy we have 
$\NN_k^{t_k}=\LG_k$ for all $k$. 
\end{cor}

The assumption $m\geq n$ in the above corollary is necessary, 
but if we require all $\NN_k$ to be prime then it can be omitted:  

\begin{cor}
Let $\G_1,\cdots,\G_n\in\CC$. 
Assume that $m\in\N$ and $\NN_1,\cdots,\NN_m$ 
are type $\mathrm{II}_1$ prime factors such that 
\[
\NN_1\ov\cdots\ov\NN_m=\LG_1\ov\cdots\ov\LG_n. 
\]
Then $m=n$ and there exist 
$t_1, t_2, ..., t_n > 0$ 
with $t_1t_2\cdots t_n = 1$ such that after permutation of indices 
and unitary conjugacy we have 
$\NN_k^{t_k}=\LG_k$ for all $k$.  
\end{cor}

Note that the above corollary shows 
the uniqueness of prime factorization 
within a certain class of type $\mathrm{II}_1$ factors, 
thus answering a question in \cite{geprime}. In 
particular, it shows that if 
$\LL\F_{r_1}\ov\cdots\ov\LL\F_{r_n}$ 
is isomorphic to 
$\LL\F_{s_1}\ov\cdots\ov\LL\F_{s_m}$, for some $2\leq r_i, s_j \leq \infty$,  
then $m=n$. This is in fact still true after tensoring by $\RR$, 
as the next corollary shows. Note in this respect 
the analogy with the 
similar results for orbit equivalence 
relations coming from 
actions of the groups $\F_{r_1} \times 
\cdots \times \F_{r_n}$ (resp. $S_\infty \times \F_{r_1} \times 
\cdots \times \F_{r_n}$) of Gaboriau \cite{ga}, and for type II$_1$ 
factors associated with actions of such groups in \cite{popab}.  

\begin{cor}
Let $\G_1,\cdots,\G_n\in\CC$.  
Assume that $m\geq n$ and $\NN_1,\cdots,\NN_m$ 
are McDuff type $\mathrm{II}_1$ factors such that
\[
\NN_1\ov\cdots\ov\NN_m=\RR\ov\LG_1\ov\cdots\ov\LG_n.
\]
Then $m=n$ and there exists a decomposition 
$\RR=\RR_1 \ov \cdots \ov \RR_n$ 
such that after permutation of indices and unitary conjugacy we have 
$\NN_k=\RR_k\ov\LG_k$ for all $k$. 
\end{cor}

\begin{cor}
Let $\G_1,\cdots,\G_n\in\CC$ and $\G=\G_1\times\cdots\times\G_n$.  
Then the fundamental group of $\LG$ is the 
product of the fundamental groups of $\LG_k$, 
i.e., $\FF(\LG) = \FF(\LG_1) \cdots \FF(\LG_n).$ 
\end{cor}

Related to the above corollary, note also that the 
outer automorphism group $\Out(\LG)$ is roughly given by 
the kernel of the map 
$\MOD \colon \Out(\LG_1^\infty) \times \cdots \times \Out(\LG_n^\infty) 
\rightarrow \R_+^*$ 
defined by 
$\MOD(\theta_1, \ldots, \theta_n) = \Mod(\theta_1) \cdots \Mod(\theta_n)$ 
($\Out(\LG)$ in fact also contains 
permutations due to possible isomorphisms between 
$\LG_k^{t_k}$, $1\leq k \leq n$). 

By using the same 
arguments, we will in fact prove generalizations 
of Theorems \ref{split1} and \ref{split2} 
to the case of inclusions of factors, as follows: 

\begin{thm}\label{emb1}
Let $\G_1,\cdots,\G_n\in\CC$, $\G=\G_1\times\cdots\times\G_n$ and 
$\NN \subset \LG$ be a subfactor.  
Assume $\NN= \NN_1\ov\cdots\ov\NN_m$ for some 
type $\mathrm{II}_1$ factors $\NN_1,\cdots,\NN_m$,  
at most one of which is injective.  
Then $m\le n$. 
If in addition $m=n$ and all $\NN_k$ are 
non-injective, then $\NN'\cap \LG$ is atomic and 
for each minimal projection $p\in \NN'\cap \LG$ 
there exist 
$t_1, t_2, \ldots, t_n > 0$, with $t_1t_2 \cdots t_n=\tau(p)$, 
and a projection $q \in \LG_1$ with $\tau(q)=\tau(p)$, such that  
after permutation of indices and unitary 
conjugacy we have  
$\NN_1^{t_1}\subset q\LG_1 q$ and 
$\NN_k^{t_k}\subset \LG_k$ 
for all $k \geq 2$. 
\end{thm}

Related to the above theorem, note that if 
$\G_1 \in \CC$ then for any $n \geq 1$ the diagonal 
embedding of $\G_1$ in $\G=\G_1 \times \cdots \times \G_1$ 
gives rise to a subfactor $\LG_1 \subset \LG$ 
with trivial relative commutant. 

\begin{thm}\label{emb2}
Let $\G_1,\cdots,\G_n\in\CC$ and  
$\NN \subset \RR\ov\LG$ be a subfactor. 
Assume that $\NN = \NN_0\ov\NN_1\ov\cdots\ov\NN_m$ 
for a type $\mathrm{II}_1$ hyperfinite factor $\NN_0$ and 
some non-injective McDuff type $\mathrm{II}_1$ factors 
$\NN_1,\cdots,\NN_m$. 
Then $m\le n$. 
Moreover, if $m=n$ and $\NN'\cap(\RR\ov\LG)=\C1$ 
then there exist mutually commuting subfactors 
$\RR_0,\RR_1,\ldots, \RR_n$ in $\RR$ 
such that 
after permutation of the indices $\{1,\ldots,n\}$ 
and unitary conjugacy we have 
$\NN_0\subset\RR_0$ and $\NN_k\subset \RR_k\ov\LG_k$ for all $k\geq 1$. 
\end{thm}

One should note that while 
the result in \cite{ozawa} 
can be viewed as a von Neumann algebra analogue 
of a theorem of Adams \cite{adams2} on the primeness 
of free ergodic m.p.\ 
actions of non-amenable hyperbolic groups,  
the above results are 
somewhat analogue to the Monod-Shalom \cite{ms} 
rigidity results for doubly ergodic actions of products 
of such hyperbolic groups. Although the 
strategies of proofs in the von Neumann 
context in \cite{ozawa} and in the present 
paper are completely different from their ergodic theory 
counterparts in \cite{adams2}, \cite{ms} 
(or for that matter \cite{ga} as well), 
it is worth mentioning 
the crucial importance for both the von Neumann algebra 
and the ergodic theory settings 
of the  amenability of the boundary action 
\cite{adams1}, \cite{ar}, \cite{ms}. Note also the similarity between the class of  
groups that 
could be dealt with in this paper and in \cite{ozawa} 
(i.e, groups in the class $\CC$ and their products), and the 
class of ``manageable'' groups in \cite{popab}, where 
however the case of lattices in $\mathrm{Sp}(1,n)$ was left unsolved 
(yet seems feasible by using \cite{cz} ``in lieu'' 
of relative Haagerup property).   

\section{Some Results on Group Algebras}\label{grp}
For a discrete group $\G$,
we denote by $\lambda$ (resp.\ $\rho$) the left (resp.\ right) 
regular representation on $\eg$ and let 
$\clg$ (resp. $\crg$) be the $C^*$-algebra generated 
by $\lambda(\G)$ (resp.\ $\rho(\G)$) in $\B(\eg)$ and 
$\LG=(\clg)''$ be its weak closure. 
The $C^*$-algebra $\clg$ (resp. the von Neumann algebra $\LG$) 
is called the reduced group $C^*$-algebra 
(resp.\ the group von Neumann algebra). 
Let $\pi_\G$ be the quotient from $\B(\eg)$ onto the 
Calkin algebra $\B(\eg)/\K(\eg)$. 
We denote by $\clrg=C^*(\lambda(\G),\rho(\G))$ 
the $C^*$-subalgebra in $\B(\eg)$ 
generated by $\clg$ and $\crg$. 

The following theorem follows from \cite{ao}, \cite{skandalis}, \cite{hg}. 

\begin{thm}\label{aoshg}
For $\G\in\CC$, the $C^*$-algebra 
$\clrg$ is exact and the $*$-homomorphism 
\[
\nu_\G\colon\clg\otimes\crg\ni\sum_{i=1}^m a_i\otimes x_i\mapsto 
\pi_\G(\sum_{i=1}^m a_ix_i)\in(\clrg+\K(\eg))/\K(\eg)
\]
is continuous w.r.t.\ the minimal tensor norm on $\clg\otimes\crg$. 
\end{thm}
\begin{proof}
The claims were already proved in Lemma 6.2.8 in \cite{hg} 
except that $\clrg$ is exact. 
So, we just sketch the proof of exactness (and continuity). 
Let $\partial\G$ be a $\G$-amenable boundary which is small at infinity
(e.g.\ the Gromov boundary for a hyperbolic group \cite{adams1}\cite{ar}). 
We identify $C(\partial\G)$ with the $C^*$-subalgebra in 
$\ell_\infty\G/c_0\G\subset\B(\eg)/\K(\eg)$. 
As it was observed in \cite{hg}, $C(\partial\G)$ commutes with 
$\crg$ in $\B(\eg)/\K(\eg)$. 
Hence, the $*$-homomorphism $\nu_\G$ extends to a $*$-homomorphism 
$\tilde{\nu}_\G$ on 
$(\G\ltimes C(\partial\G))\otimes(\G\ltimes C(\partial\G))^{\mathrm{op}}$. 
Since $\G\ltimes C(\partial\G)$ is nuclear, 
$\tilde{\nu}_\G$ is continuous w.r.t.\ the minimal tensor norm. 
Moreover, $\pi_\G(\clrg)$ is contained 
in a nuclear $C^*$-subalgebra of $\B(\eg)/\K(\eg)$. 
Thus, by the Choi-Effros lifting theorem \cite{ce}, 
$\pi_\G|_{\clrg}$ has a unital completely positive lifting. 
This implies the exactness of $\clrg$. 
\end{proof}

We note that since we assumed $\G$ to be non-amenable and ICC, 
the ideal $\K(\eg)$ is actually contained in $\clrg$ 
(otherwise $\clrg\cong\clg\otimes_{\min}\crg$ and hence $\G$ is amenable) 
and the $*$-homomorphism $\nu_\G$ is injective. 

 From now on, 
let $\G_0$ be either the trivial group or an amenable ICC group. 
We note that $\clrg_0\cong\clg_0\otimes_{\min}\crg_0$. 
Let $\G_1\ldots,\G_n\in\CC$ 
and $\G=\G_0\times\G_1\times\cdots\times\G_n$. 
We note that $\clg$ is locally reflexive (and exact). 
For each $k$, let $\G_k'$ be the kernel of 
the projection of $\G$ onto $\G_k$, i.e., 
$\G=\G_k'\times\G_k$. 
There is an obvious identification 
\[
\clrg=\clrg_0\om\clrg_1\om\cdots\om\clrg_n
\subset\B(\eg)=\B(\bigotimes_k\eg_k). 
\]
Since each $\clrg_k$ is exact by Theorem \ref{aoshg}, 
the kernel of the $*$-homomorphism
\begin{align*}
& \clrg_0\om\clrg_1\om\cdots\om\clrg_n\\
& \qquad \to\clrg_0\om(\clrg_1/\K(\eg_1))\om\cdots\om(\clrg_n/\K(\eg_n))
\end{align*}
is 
\[
J=\clrg_1'\om\K(\eg_1)+\cdots+\clrg_n'\om\K(\eg_n)\subset\clrg.
\]
We denote by $\pi_J\colon\clrg\to\clrg/J$ the quotient map. 
The following lemma summerizes the above arguments. 
\begin{lem}
For $\G=\G_0\times\G_1\times\cdots\times\G_n$ as above, 
the $*$-homomorphism 
\[
\nu_\G\colon\clg\otimes\crg\ni\sum_{i=1}^m a_i\otimes x_i
\mapsto\pi_J(\sum_{i=1}^m a_ix_i)\in\clrg/J
\]
is continuous w.r.t.\ the minimal tensor norm on $\clg\otimes\crg$. 
\end{lem}

Note that by Connes' theorem \cite{connes}, 
we may find an amenable subgroup $\Delta$ 
in the unitary group $\U(\RR)$ of the injective factor $\RR$ 
which generates $\RR$ as a von Neumann algebra. 
Thus, whenever $\RR$ is represented 
as a von Neumann algebra on a Hilbert space, $\RR\subset \B(\hh)$, 
by averaging $\Ad u$ over $u\in\Delta$ along some F{\o}lner sequence, 
we obtain a conditional expectation $\Psi_{\RR'}$ 
from $\B(\hh)$ onto $\RR'$ such that 
\[
\Psi_{\RR'}(a)\in\conv\{uau^* : u\in\U(\RR)\}.
\]
Such a conditional expectation is said to be proper. 
\begin{prop}\label{pos}
Let $\G=\G_0\times\G_1\times\cdots\times\G_n$ be as above 
and let $\RR\subset\LG$ be an injective von Neumann subalgebra 
with a proper conditional expectation 
$\Psi_{\RR'}$ from $\B(\eg)$ onto $\RR'$ in $\B(\eg)$. 
If the relative commutant $\RR'\cap\LG$ is non-injective, 
then there exists $k\in\{1,\ldots,n\}$ such that  
\[
\Psi_{\RR'}(\C1_{\eg_k'}\otimes\K(\eg_k))\neq\{0\}.
\]
\end{prop}
\begin{proof}
Let $\NN=\RR'\cap\LG$ and let $\e_{\NN}$ be 
the trace preserving conditional expectation from $\LG$ onto $\NN$. 
Then, we have $\Psi_{\RR'}|_{\LG}=\e_{\NN}$ by uniqueness. 

For a proof by contradiction, suppose that 
$\Psi_{\RR'}(\C1_{\eg_k'}\otimes\K(\eg_k))=\{0\}$ for all $k$. 
It follows that $J\subset\ker\Psi_{\RR'}$ and hence 
$\Psi_{\RR'}|_{\clrg}=\tilde{\Psi}_{\RR'}\circ\pi_J$ 
for some unital completely positive map 
$\tilde{\Psi}_{\RR'}\colon\clrg/J\to\B(\eg)$. 
It follows that the unital completely positive map 
$\Phi_{\NN}
=\tilde{\Psi}_{\RR'}\circ\nu_\G\colon\clg\otimes_{\min}\crg\to\B(\eg)$ 
is continuous and we have 
\[
\Phi_{\NN}(\sum_{i=1}^m a_i\otimes x_i)
=\sum_{i=1}^m \e_{\NN}(a_i)x_i\in\B(\eg).
\]
for $\sum_{i=1}^m a_i\otimes x_i\in \clg\otimes\crg$. 
However, by Lemma~5 in \cite{ozawa}, this implies 
the injectivity of $\NN$ which contradicts the assumption.  
\end{proof}
\section{Some Results on Subfactors and Unitary Conjugacy}
\begin{prop}\label{split}
Let $\MM=\MM_1 \ov \MM_2$ 
and $\NN \subset \MM$ be type $\mathrm{II}_1$ factors. 
Assume there exists 
$x\in \MM_1\otimes\K(L^2(\MM_2))\subset\B(L^2(\MM))$ 
such that 
\[\NN'\cap\conv\{uxu^* : u\in \U(\NN)\}\neq\{0\}.\] 
Then there exist a decomposition 
$\MM=\MM_1^t \ov \MM_2^{1/t}$ for some $t >0$, 
and a non-zero partial isometry $u \in \MM$ such that 
$uu^*\in \Proj(\MM_2^{1/t})$, $u^*u \in \Proj(\NN'\cap \MM)$ 
and $u\NN u^* \subset \MM_1^t uu^*$. 
Moreover, if $\NN'\cap\MM$  is a factor, 
then $u$ can be taken a unitary element. 
\end{prop}
\begin{proof}
By Proposition 1.3.2 in \cite{popab}, 
$\NN'\cap\conv\{uxu^* : u\in \U(\NN)\}$ 
consists in just one element $y$, 
which is contained in $\K(L^2(M_1))\otimes_{\min}\MM_2$. 
Thus, if $y \neq 0$ and we denote by $e$ the spectral 
projection of $|y|$ corresponding to the interval 
$[\|y\|/2, \|y\|]$ then $e$ is a finite non-zero projection in 
$\MM_1 \ov \B(L^2(\MM_2))$, which commutes with $\NN$. 

Thus, the Hilbert space $\hh=eL^2(\MM)$ is an $\NN-\MM_1$ 
Hilbert bimodule and $\dim\hh_{\MM_1} < \infty$. 
Thus, if we take $p\in \Proj(\MM_1)$, $q_0\in \Proj(\NN)$ with 
$\tau(p)= \min\{s, 1\}$, $\tau(q_0) =\min\{1/s, 1\}$, 
where $s=\dim\hh_{M_1}$, then 
$\hh_0 = q_0\hh p$ is a $q_0\NN q_0-p\MM_1p$ Hilbert bimodule 
and $\dim(\hh_0)_{p\MM_1p}=1$. 

Thus, there exists a vector $\xi_0 \in q_0L^2(\MM)p$ such that 
$\hh_0 = L^2(\xi_0 p\MM_1p)=\overline{\xi_0 p\MM_1p}$. 
Since $q_0\NN q_0\hh_0 p\MM_1p = \hh_0$, this implies 
$q_0\NN q_0 \xi_0 \subset L^2(\xi_0 p\MM_1p)$. 

By regarding $\xi_0$ as a square summable operator affiliated 
with $\MM$ and using the fact that $x\in p\MM_1p$, 
$\xi_0 x=0$ implies $x=0$, it follows that the summable operator 
$\e_{\MM_1}(\xi_0^*\xi_0)\in L^1(\MM_1)$ 
affiliated with $\MM_1$ has support $p$. 
Thus, $\xi = \xi_0 \e_{\MM_1}((\xi_0^*\xi_0)^{-1/2})$ is 
a square summable operator affiliated with $\MM_1$ such that 
$\e_{\MM_1}(\xi^*\xi)=p$. 
In other words, as an unbounded operator acting on $L^2(\MM)$, 
the restriction $\xi|_{L^2(p\MM_1p)}$ is a partial isometry. 

We then have $L^2(\xi p\MM_1p) = \xi L^2(p\MM_1p)$. 
Let $q\in q_0\NN q_0$ be the minimal projection in $q_0\NN q_0$ 
with the property that $(q_0-q)\xi=0$. 
We denote $\theta(x)=\e_{\MM_1}(\xi^*x\xi)$, $x\in q\NN q$, 
and note that $\theta$ is a unital, normal, faithful, 
completely positive map from $q\NN q$ into $p\MM_1p$. 
Moreover, the same argument as in Lemma 5 of \cite{popaf} 
shows that in fact $\theta$ is a unital $*$-isomorphism 
of $q\NN q$ into $p\MM_1p$. 

Thus, since $\xi$ is an intertwiner between $\id|_{qNq}$ and  
the isomorphism $\theta$, if $v\in\MM$ is the 
partial isometry in the polar decomposition of $\xi$,
then $v^*v \in (q\NN q)'\cap q\MM q$, 
$vv^* \in \theta(q\NN q)'\cap p\MM p$ 
and $vx=\theta(x)v$ for all $x\in q\NN q$. 

Noticing that if $\BB\subset p\MM_1p$ is a von Neumann subalgebra 
then $\BB'\cap p\MM p=(\BB'\cap p\MM_1p) \vee \MM_2$ and that 
\[
\ZZ(\BB'\cap p\MM p) = \ZZ(\BB'\cap p\MM_1p) 
= \ZZ(\ZZ(\BB'\cap p\MM_1p) \vee \MM_2).
\]
It follows that the projection $vv^* \in \theta(q\NN q)'\cap p\MM p$ 
is equivalent in $\theta(q\NN q)'\cap p\MM p$ with a projection of 
same central trace in $\ZZ \vee \MM_2$, where 
$\ZZ = \ZZ(\theta(q\NN q)'\cap p\MM_1p) \subset p\MM_1p$. 
Thus, we may assume the partial isometry $v$ 
is so that $vv^* \in \ZZ \vee \MM_2$. 
Then there exist non-zero projections 
$p_1 \in \ZZ \subset p\MM_1p$, 
$p_2\in \MM_2$ such that the central trace 
$\tau(p_2)p_1$ of $p_1p_2$ in $\ZZ \vee \MM_2$ 
is less than or equal to the central trace of $vv^*$, 
i.e., there exists $w\in \ZZ \vee \MM_2$ 
such that $ww^* = p_1p_2$, $w^*w \leq vv^*$. 

But then the partial isomtery $v_0=wv$ satisfies 
$v_0v_0^* =p_1p_2$ and 
\[
v_0 x = wv x= w \theta(x)v = \theta (x)wv =\theta(x)v_0,  
\]
in which we can now regard $\theta$ as a 
unital $*$-isomorphism of $q\NN q$ into $p_1\MM_1p_1$. 
Also, we still have $v_0^*v_0\in (q\NN q)'\cap q\MM q$, 
i.e., there exists a projection $q'\in\NN'\cap\MM$ 
such that $v_0^*v_0=qq'$. 
Moreover, by restricting $\sigma$ and $v_0$ 
to a projection of appropriate trace in $\Proj(q\NN q)$, 
we may assume $\tau(q) = 1/n,$ for some integer $n \geq 1$.  

Let now $t = \tau(p_1)/\tau(q)$ and decompose $\MM$ 
as $\MM_1^t \ov \MM_2^{1/t}$. 
Since $t \geq \tau(p_1)$, it follows that 
we may assume $p_1\MM_1p_1 \subset \MM^t$ and  
$p'=p_2 \in \MM_2^{1/t}$.  
Moreover, by the choice of $t$ the normalized 
traces on $\NN$ and $\MM^t$ satisfy 
$\tau_{\NN}(q) = \tau_{\MM_1^t}(p_1)=1/n$. 
Let then $w_1, w_2, ..., w_n \in\NN$, 
$u_1, u_2, ..., u_n \in \MM_1^t$ be partial 
isometries such that $w_iw_i^*=q$, $u_i^*u_i =p_1$, for any $i$, 
and $\sum_i w_i^*w_i = 1_N$, $\sum_i u_iu_i^* = 1_{\MM_1^t}$. 

Thus, if we let  $u=\sum_i u_iv_0w_i$ 
then $u^*u=q' \in \NN'\cap \MM$, $uu^*=p' \in \MM_2^{1/t}$ 
and $u\NN u^* \subset \MM_1^tp'$. 

Similarly, if $\NN'\cap\MM$ is a factor then there 
exist partial isometries $w'_1, w'_2,\ldots, w'_m \in \NN'\cap \MM$, 
$u'_1, u'_2, \ldots, u'_m \in \MM_2^{1/t}$ 
such that $w'_j{w'_j}^*\leq q'$, 
${u'_i}^*u'_i \leq p'$ for any $j$ 
and $\sum_i {w'_i}^*w'_i = 1_{\NN'\cap \MM}$, 
$\sum_i u'_i{u'_i}^* = 1_{\MM_2^{1/t}}$. 
Then the unitary element 
$u'=\sum_j u'_j u w'_j$ will satisfy $u'\NN{u'}^* \subset \MM^t$. 
\end{proof}

The following result is an easy consequence of 
A.1.1 and A.1.2 in \cite{popac}, 
but we give full details for completeness. 

\begin{prop}\label{injsf}
Let $\NN \subset \MM$ be type $\mathrm{II}_1$ factors 
with separable preduals. 
Then there exists a hyperfinite subfactor $\RR \subset \NN$ 
such that $\RR'\cap\MM = \NN'\cap\MM$. 
\end{prop}
\begin{proof}
Let $\{x_n\}_n$ be $\|\quad\|_2$-dense 
in the unit ball of $\MM$. 
We construct recursively mutually commuting 
finite dimensional subfactors $\NN_k \subset \NN$, $k=1, 2,\ldots,$ 
such that if we put 
$\MM_n = \bigvee_{1\leq k \leq n} \NN_k$ 
then 
\[
\|\e_{\MM_n'\cap\MM}(x_j)-\e_{\NN'\cap\MM}(x_j)\|_2 
\leq 2^{-n},\ \forall j \leq n. 
\]
 
Assume we have constructed $\NN_1, \NN_2, ..., \NN_n$. 
Let $y_j = \e_{\NN_n'\cap\MM}(x_j)$, $j \leq n+1$. 
By A.1.1 in \cite{popac}, given any $\varepsilon > 0$ 
there exists an abelian finite dimensional 
subalgebra $\A_0 \subset \MM_n'\cap \NN$ such that 
if we denote $\QQ=\MM_n'\cap\NN$, $\PP=\MM_n'\cap\MM$ 
and $y'_j = y_j - \e_{\QQ \vee \QQ'\cap \PP}(y_j)$, then 
$\|\e_{\A_0'\cap \PP}(y'_j)\|_2 \leq \varepsilon$ 
for all $j \leq n+1.$ 
Moreover, $\A_0$ can be taken to have 
minimal projections of equal trace. 
Thus, if we let $\NN_0 \subset \QQ$ be a finite dimensional 
factor such that $\A_0 \subset \NN_0$ is its diagonal algebra 
then we still have 
$\|\e_{\NN_0'\cap \PP}(y'_j)\|_2 \leq \varepsilon$ 
for all $j \leq n+1.$ 

Denote 
$y''_j = \e_{\QQ \vee \QQ'\cap \PP}(y_j)$ 
and let $y_{j,k} \in \QQ$, $y'_{j,k}\in \QQ'\cap \PP$, 
$1\leq k \leq n_j,\ 1\leq j \leq n+1$ 
be such that 
$\|y''_j-\sum_k y_{j,k}y'_{j,k}\|_2 \leq \varepsilon$, for all $j\leq n+1$. 
By the first part of proof of A.1.2 in \cite{popac}, 
there exists a finite dimensional factor $\MM_0 \subset \NN_0'\cap \QQ$ 
such that $\|\e_{\MM_0'\cap \QQ}(y_{j,k})-\tau(y_{j,k})1\|_2 \leq \delta$, 
where $\delta = \varepsilon/(\sum_j n_j)$. 
It follows that 
\begin{align*}
\|\e_{\MM_0'\cap \PP}(y''_j) - \e_{\QQ'\cap\PP}(y''_j)\|_2 
&\leq 
\|\e_{\MM_0'\cap \PP}(y''_j) - \sum_k \tau(y_{j,k}) y'_{j,k}\|_2\\ 
&\leq \varepsilon + n_j\delta \leq 2 \varepsilon,\ \forall j\leq n+1. 
\end{align*}
Altogether, if we denote 
$\NN_{n+1} = \NN_0 \vee \MM_0 \subset \QQ = \MM_n'\cap \NN$ 
and put $\MM_{n+1} = \MM_n \vee \NN_{n+1}$ then 
\begin{align*}
\|\e_{\MM_{n+1}'\cap \MM}(x_j) - & \e_{\NN'\cap \MM}(x_j)\|_2 
= \|\e_{\MM_{n+1}'\cap \MM}(y_j) - \e_{\NN'\cap \MM}(y_j)\|_2\\
&\leq \|\e_{\MM_{n+1}'\cap \MM}(y'_j)\|_2 
 + \|\e_{\MM_{n+1}'\cap \MM}(y''_j) - \e_{\NN'\cap \MM}(y''_j)\|_2\\
&\leq \|\e_{\NN_0'\cap \PP}(y'_j)\|_2 
 + \|\e_{\MM_0'\cap \PP}(y''_j) - \e_{\QQ'\cap \PP}(y''_j)\|_2\\
&\leq 3\varepsilon,\ \forall j \leq n+1. 
\end{align*}
Thus, if we take $\varepsilon = 2^{-n-1}/3$ then $\MM_{n+1}$ 
satisfies the required conditions. 

Finally, if we put $\RR = \overline{\bigcup_n \MM_n}^w$ then 
\[
\e_{\RR'\cap \MM}(x_j) = \lim_n \e_{\MM_n'\cap \MM}(x_j) 
= \e_{\NN'\cap \MM}(x_j),\ \forall j \geq 1
\]
and thus, by density, 
$\e_{\RR'\cap \MM}(x) = \e_{\NN'\cap \MM}(x)$, for all $x\in \MM$. 
\end{proof}

We end this section with a technical result 
on geometry of projections, needed in the proof 
of Theorem 6. 

\begin{lem}\label{proj}
Let $\NN=\NN_1\ov \NN_2 \subset \MM$ be an 
inclusion of type $\mathrm{II}_1$ factors and 
denote $\BB = \NN'\cap \MM$. If $p \in \NN_2'\cap \MM$ 
is a non-zero projection then there exist non-zero projections 
$q \in \BB \vee \NN_1, q_0 \in \BB, q_1 \in \NN_1$ 
such that $p$ is equivalent to $q$ in 
$\NN_2'\cap M$ and $p$ majorises $q_0q_1$ in 
$\NN_2'\cap M$. 
\end{lem}
\begin{proof} We clearly have 
$\ZZ(\NN_2'\cap \MM) \subset \ZZ(\BB \vee \NN_1)$. 
Let $z$ be the central trace of $p$ in $\NN_2'\cap \MM$. 
Since $\BB \vee \NN_1$ is of type II$_1$, $0\leq z \leq 1$ 
and $z \in \ZZ(\BB \vee \NN_1)$, 
it follows that there exists a projection $q \in 
\BB \vee \NN_1$ with central trace in $\BB \vee \NN_1$ 
equal to $z$. Thus, $p$ and $q$ are 
equivalent in $\NN_2'\cap \MM$. 

Moreover, 
since there exists a non-zero 
projection $q_0 \in \ZZ(B)=\ZZ(B \vee \NN_1)$ 
and $c > 0$ such that  $cq_0 \leq z$, 
it follows that if $q_1 \in \NN_1$ 
is a projection of trace $c$ then $q_0q_1$ 
has central trace $cq_0\leq z$ in $\ZZ(B \vee \NN_1)$, 
thus being majorised in $B \vee \NN_1$ 
by $q$ and in $\NN_2'\cap \MM$ by $p$ as well. 
\end{proof}
\section{Proof of the Main Theorem}
The following proposition is a combination of 
Propositions \ref{pos}, \ref{split} and \ref{injsf}. 

\begin{prop}\label{prp}
Let $\G_0$ be either a trivial group or an amenable ICC group, $\G_1,\cdots,\G_n\in\CC$ and 
$\G=\G_0\times\G_1\times\cdots\times\G_n$. 
Assume that $\NN\subset\LG$ is a type $\mathrm{II}_1$ factor 
such that $\NN'\cap\LG$ is non-injective. 
Then, there exist $k\in\{1,\ldots,n\}$, 
$t>0$ and a non-zero partial isometry $u\in\LG$ 
with $e=uu^*\in \Proj((\LG_k)^{1/t})$, 
$u^*u \in \Proj(\NN'\cap\LG)$ 
such that 
\[
u\NN u^*\subset(\LG_k')^t\otimes\C e
\subset(\LG_k')^t\ov e(\LG_k)^{1/t}e
=e(\LG)e.
\]
Moreover, if $\NN'\cap\LG$  is a factor, 
then $u$ can be taken a unitary element. 
\end{prop}
\begin{proof}
We denote $\MM=\LG$. 
By Proposition~\ref{injsf}, 
there exists a hyperfinite subfactor $\RR\subset\NN$ 
such that $\RR'\cap\MM=\NN'\cap\MM$. 
By Proposition~\ref{pos}, 
there exist $k\in\{1,\ldots,n\}$ and 
$x\in\C1_{\eg_k'}\otimes\K(\eg_k)$ 
such that 
\[
\RR'\cap\conv\{ uxu^* : u\in\U(\RR)\}\neq\{0\}.
\]
We denote $\MM_1=\LG_k'$ and $\MM_2=\LG_k$. 
Applying Proposition~\ref{split} to $\RR\subset\MM$, 
we obtain $t>0$ and and a partial isometry $u\in\LG$ 
with $e=uu^*\in \Proj((\LG_k)^{1/t})$, 
$u^*u \in \Proj(\NN'\cap\LG)$ 
such that 
$u\RR u^*\subset\MM_1^t e$. 
In case where $\NN'\cap\MM$ is a factor, 
$u$ can be taken a unitary element. 
It follows that 
\begin{align*}
u\NN u^* &\subset (u(\NN'\cap\MM)u^*)'\cap(u\MM u^*)\\
&= ((u\RR u^*)'\cap(u\MM u^*))'\cap(u\MM u^*)\\
&\subset ((\MM_1^te)'\cap(e\MM e))'\cap(e\MM e)
= \MM_1^te. 
\end{align*}
This completes the proof.
\end{proof}

Now we are in a position 
to prove Theorems \ref{split1} and \ref{split2}. 
\begin{proof}[Proof of Theorems \ref{split1} and \ref{split2}]
We proceed by induction over $n$. 
The case where $n=0$ is trivial. 
Let $\G_0$ be the trivial group (resp.\ an amenable ICC group), 
$n\geq 1$, $\G_1,\ldots,\G_n\in\CC$ and let 
$\G=\G_0\times\G_1\times\cdots\times\G_n$. 
Since $\MM=\LG$ is non-injective, 
we may assume, say, $\NN_2$ is non-injective. 
It follows from Proposition~\ref{prp} 
that there exist $t>0$, a unitary element $u\in\LG$ 
and $k\in\{1,\ldots,n\}$ such that 
if we denote $\MM_1=\LG_k'$ and $\MM_2=\LG_k$ then 
$u \NN_1 u^*\subset \MM_1^t$. 
Thus, with $\NN_{2,1}=\NN_1'\cap u^* \MM_1^t u$, 
we have 
\[
\NN_2 = \NN_1'\cap\MM = u^*\big( \NN_{2,1}\ov\MM_2^{1/t} \big)u
\subset u^*\big( \MM_1^t\ov\MM_2^{1/t} \big)u =\MM
\]
and $\NN_1\vee\NN_{2,1}=u^*\MM_1^tu$. 
Therefore, by the induction hypothesis, 
there exist $s>0$, 
a partition $I_1\cup I_{2,1}=\{1,\ldots,n\}\setminus\{k\}$ 
and $\RR_0=\RR_1=\C1$ 
(resp.\ a decomposition $\RR = \RR_1 \ov \RR_2$)
such that 
$\NN_1 = (\RR_1\ov\overline{\bigotimes}_{k\in I_1}\LG_k)^s$ 
and 
$\NN_{2,1}= (\RR_2\ov\overline{\bigotimes}_{k\in I_{2,1}}\LG_k)^{t/s}$ 
modulo conjugacy by a unitary element in $(\LG_k')^t$. 
We complete the proof by setting $I_2=\{k\}\cup I_{2,1}$. 
\end{proof}

An adaptation of the above proof 
shows that if $\RR\ov\NN=\RR\ov\LG$ 
for $\G\in\CC$ and some 
non-McDuff type $\mathrm{II}_1$ factor $\NN$, 
then there exists $t>0$ such that 
$\NN^t=\LG$ in $\RR\ov\LG$ modulo unitary conjugacy. 

\begin{prop}\label{num}
Let $\G_0$ be the trivial group (resp.\ an amenable ICC group), 
and $\G_1,\cdots,\G_n\in\CC$. 
Let $\NN_1$ be a diffuse (resp.\ non-injective) von Neumann algebra, 
$\NN_2,\cdots,\NN_m$ be non-injective type $\mathrm{II}_1$ factors. 
Assume that 
\[
\NN_1\ov\cdots\ov\NN_m\subset\LG_0\ov\LG_1\ov\cdots\ov\LG_n.
\]
Then $m\le n$. 
\end{prop}
\begin{proof}
We proceed by induction over $n$. 
The case where $n=0$ is trivial. 
We first deal with the case where $\G_0$ is non-trivial 
and $\NN_1$ is non-injective. 
In this case, by Proposition~\ref{prp} 
there exist $k\in\{1,\ldots,n\}$ and $s>0$ such that 
$(\NN_2\ov\cdots\ov\NN_m)^s$ 
is unitary conjugate to a subfactor of $\LG_k'$, 
and hence we obtain $m-1\le n-1$ by the induction hypothesis. 

Now, let $n\geq 1$, $\G=\G_1\times\cdots\times\G_n$ 
and let $\NN_1,\ldots,\NN_m$ be as in the statement. 
For a proof by contradiction, suppose $m>n$ and 
let $\PP=\NN_2\ov\cdots\ov\NN_{m-1}$ (or just $\C1$ if $m=2$). 
By Proposition~\ref{prp}, there exist 
$k\in\{1,\ldots,n\}$, 
$t>0$ and a partial isometry $u\in\LG$ 
with $e=uu^*\in \Proj((\LG_k)^{1/t})$, 
$u^*u \in \Proj(\PP'\cap\LG)$ 
such that 
\[
u\PP u^*\subset(\LG_k')^t\otimes\C e
\subset(\LG_k')^t\ov e(\LG_k)^{1/t}e
=e(\LG)e.
\]
By Lemma \ref{proj}, we may assume that 
$uu^*=q_0q_1$ with $q_0\in\BB=(\NN_2\ov\cdots\NN_m)'\cap\LG$
and $q_1\in\NN_m$. 
Since $\PP$ is unitary conjugate to a subfactor in $(\LG_k')^t$, 
it follows from the induction hypothesis and the assumption $m>n$ 
that $m-2=n-1$. 
Thus, by the induction hypothesis again, 
$\A=(u\PP u^*)'\cap(\LG_k')^t$ 
cannot be diffuse. 
Hence for a minimal projection $f\in \A$, 
we have 
\[
fu(q_0\BB q_0\vee q_1\NN_mq_1)u^*f 
\subset \C f \otimes e(\LG_k)^{1/t}e 
\subset f(\LG_k')^tf \ov e(\LG_k)^{1/t}e.
\]
But the first inclusion contradicts the solidity of $\LG_k$ 
proved in \cite{ozawa}. 
\end{proof}

\begin{proof}[Proof of Theorems \ref{emb1} and \ref{emb2}] 
In each statement, the first part follows from Proposition~\ref{num}. 

Now we prove the second part by induction over $m=n$. 
The case where $n=0$ is trivial. 
Let $\G_0$ be the trivial group 
(resp.\ an amenable ICC group) and $\G_1,\cdots,\G_n\in\CC$. 
Assume that $\NN_0=\C1$ (resp.\ $\NN_0$ is a hyperfinite factor) 
and $\NN_1,\ldots,\NN_n$ are non-injective (resp.\ non-injective McDuff) 
factors such that 
\[
\NN=\NN_0\ov\NN_1\ov\cdots\ov\NN_n
\subset\LG=\LG_0\ov\LG_1\ov\cdots\ov\LG_n. 
\]
Then, by Proposition~\ref{num}, 
$\NN'\cap\MM$ is atomic (resp.\ injective). 
To prove the rest, we may assume that $\NN'\cap\MM=\C1$. 
Let $\PP=\NN_0\ov\NN_1\ov\cdots\ov\NN_{n-1}$ 
(or $\PP=\NN_0$ if $n=1$). 
By Proposition~\ref{prp},
there exist $k\in\{1,\ldots,n\}$, $t>0$ 
and a unitary element $u$ in $\LG$ 
such that 
\[
u \PP u^* \subset (\LG_k')^t\otimes\C 1
\subset(\LG_k')^t\ov(\LG_k)^{1/t}=\LG.
\]
Since $\NN'\cap\LG$ is a factor, so is 
$\QQ=(u\PP u^*)'\cap(\LG_k')^t$, 
but then $\QQ$ has to be 
finite dimensional (resp.\ hyperfinite) 
by Proposition~\ref{num}. 
We have 
\[
u \NN_n u^* \subset (u\PP u^*)'\cap \LG 
= \QQ\ov(\LG_k)^{1/t}.
\]
Therefore, we have $\NN_n\subset(\LG_k)^{d/t}$ with 
$\dim\QQ=d^2$ (resp.\ $\NN_n\subset\QQ^{1/t}\ov\LG_k$).
By the induction hypothesis, 
this finishes the proof of Theorem \ref{emb1}. 
For Theorem~\ref{emb2}, one has to apply the induction hypothesis 
to the subfactor 
\[
(\QQ^{1/t}\vee\NN_0)\ov\NN_1\ov\cdots\ov\NN_{n-1}\subset \LG_k'
\]
and obtain $\QQ^{1/t}\subset\LG_0=\RR$. 
\end{proof}
\noindent\textbf{Acknowledgment.}
This research was carried out while the first named 
author was visiting 
the University of California at Los Angeles under the support of 
the Japanese Society for the Promotion of Science
Postdoctoral Fellowships for Research Abroad.
He gratefully acknowledges the kind hospitality of UCLA. 

\end{document}